\documentclass[12pt,a4paper]{article}
\usepackage[english]{babel}
\usepackage{amssymb}
\usepackage{amsmath}
\usepackage[dvips]{graphicx}
\def\CC{\mathbb{C}}

\def\ZZ{\mathbb{Z}}

\def\A{{\mathcal A}}

\def\D{{\mathcal D}}
\def\E{{\mathcal E}}

\def\M{{\mathcal M}}

\def\R{{\mathcal R}}
\def\K{{\mathcal K}}

\def\g{\gamma}
\def\e{\varepsilon}
\def\O{{\mathcal O}}
\def\SS{{\mathcal S}}

\def\a{\alpha}
\def\b{\beta}
\def\d{\delta}

\def\o{\omega}
\def\cs{{\check s}}
\def\cM{{\check M}}

\def\l{\lambda}

\def\f{\varphi}

\def\p{\partial}

\def\Mat{\operatorname{Mat}}
\def\Pf{\operatorname{Pf}}
\def\SL{\operatorname{SL}}
\def\GL{\operatorname{GL}}
\def\PGL{\operatorname{PGL}}
\def\Gr{\operatorname{Gr}}
\def\Sk{\operatorname{Sk}}
\def\Sym{\operatorname{Sym}}
\def\Sq{\operatorname{Sq}}
\def\tr{\operatorname{tr}}

\def\diag{\operatorname{diag}}
\def\min{\operatorname{min}}
\newtheorem{theorem}{Theorem}[section] 
\newtheorem{corollary}[theorem]{Corollary}

\newtheorem{example}[theorem]{Example} 
\newtheorem{prop}[theorem]{Proposition}
\newcommand{\bprop}{\begin{prop}}
\newcommand{\eprop}{\end{prop}}
\newtheorem{definition}[theorem]{Definition}
\newtheorem{rem}[theorem]{Remark}
\newcommand{\brem}{\begin{rem}}
\newcommand{\erem}{\end{rem}}
\newtheorem{rems}[theorem]{Remarks}
\newcommand{\brems}{\begin{rems}}
\newcommand{\erems}{\end{rems}}

\begin{document}

\title{{\bf Families of skew-symmetric matrices \\  of even size}}

\author{J.W. Bruce, V.V.  Goryunov and G.J. Haslinger}
\date{}
\maketitle

\begin{abstract}
The main result of the paper is a classification of singularities of complex skew-symmetric matrix families of even size which are simple under a natural equivalence relation. The classification is obtained by appropriate suspensions of simple families of arbitrary square matrices classified by Bruce and Tari. The suspensions we are using are based on the embedding of a space of arbitrary square matrices into the space of skew-symmetric matrices of twice the size. We also show that similar relations embed Bruce's classification of simple symmetric matrix families into the Bruce-Tari simple list. Our constructions introduce a unified approach to all three simple classifications.    
\end{abstract}

\ 

Matrices depending on parameters are traditional objects of numerous mathematical areas ranging from the most applied to the most abstract. So, an understanding of the local behaviour of matrix families has always been a natural and important problem. On the other hand, a starting point of any local understanding is a classification of singularities which are simple (that is, have no moduli) with respect relevant equivalences.   

This paper is a 
continuation
(promised quite a while ago)
of work on the classification of simple singularities of three types of families of square matrices --- arbitrary, symmetric and skew-symmetric --- under their natural equivalences (see Section \ref{Seq} for details). We denote the equivalences and family types respectively $\Sq$, $\Sym$ and $\Sk$. 

The work was started by the first author with his paper \cite{B} on the classification of $\Sym$-simple families. It was followed by his paper with Tari \cite{BT} in which $\Sq$-simple singularities were classified, and by a thesis of the third author \cite {H} where a list of two-parameter simple $4 \times 4$ skew-symmetric matrix families was obtained. 

This paper gives a complete classification of complex $\Sk$-simple families of arbitrary even size.  This classification turns out to be very closely related to the $\Sq$-classification from \cite{BT}. The relation is based on the embeddings $M \mapsto \left( \begin{array}{cc} 0 & M \\ -M^\top & 0 \end{array} \right)$ of spaces of square matrices into the spaces of skew-symmetric matrices of double the size. We introduce three ways --- we call them {\em trivial, full} and {\em  intermediate suspensions} --- to produce a skew-symmetric family from an arbitrary square one, and it emerges from the calculations that these are sufficient to obtain our list from that in \cite{BT}. Our suspensions work like stabilisations of functions by addition of squares of new variables. The trivial suspension approach already appeared in \cite{H}. 

On the linear level, the trivial suspensions relate pencils of arbitrary $k \times k$ matrices and of skew-symmetric $2k \times 2k$ matrices (see \cite{T}), while the duals of such pencils are related by the full suspensions. However, existence of extensions of these correspondences from the 1-jets to the entire simple classifications was quite unexpected.  

Comparison of the simple lists from \cite{B} and \cite{BT} shows that the standard embeddings of spaces of symmetric matrices into the spaces of arbitrary square matrices of the same size provide two further suspension types which embed the $\Sym$-simple classification into the $\Sq$ one. Thus, the set of all five types of suspensions provides a unified point of view on all three simple classifications. 

\

The paper is organised as follows.

Section \ref{Seq} recalls the equivalences of arbitrary square, symmetric and skew-symmetric matrix families with which we are working in this paper.

Section \ref{Ssq} reminds the reader of the classification of simple square matrix families obtained in \cite{BT}.

In Section \ref{Ssus} we introduce our three types of the $\Sq$-to-$\Sk$ suspensions.
We establish relations between the classes of the initial and suspended singularities. In Section \ref{Smain}, we formulate our main result, Theorem \ref{Tmain}, which gives a complete classification of simple complex even-size skew-symmetric matrix families in terms of the suspensions of simple arbitrary square matrix singularities.

In Section \ref{Sdim} we take our first step towards a proof of Theorem \ref{Tmain} and derive constraints on the size of skew-symmetric matrices and on the number of parameters in a family that may allow existence of simple singularities.
All the conditions obtained  are collected in Corollary \ref{Cdim}.  

Section \ref{Sred} completes our proof of Theorem \ref{Tmain}. There we reduce the classification of simple skew-symmetric map germs of the fixed matrix corank 
in each of the cases allowed by Corollary \ref{Cdim} to the suspensions of half-the-size simple square matrix classification according to the claims of Theorem \ref{Tmain}. In the final subsection we check that none of the elements of these skew-simple  classifications is adjacent to non-simple classes of a smaller matrix corank.
At this step the criterion that simplicity of a matrix family is equivalent to the simplicity of its determinantal or Pfaffian function, in the dimensions where these functions have isolated singularities, turns out to be very useful. 

Finally, in Section \ref{Ssymtosq}, we comment  on the suspension embedding of the $\Sym$-simple classification into the $\Sq$-simple, and bring all three simple classifications together. 

We work throughout over the complex numbers. Lists of simple matrix singularities over the reals can be derived using the same approach.

\ 

\section{Matrix equivalences} \label{Seq}

A {\em matrix family\/} 
in this paper will be  a holomorphic mapping $M: U \to \Mat_n,$ where $U$ is an open domain in  $\CC^s,$ and $\Mat_n$ is one of three spaces of complex $n \times n$ matrices:
\begin{itemize}
\item[] $\Sq_n,$ arbitrary square matrices,
\item[] $\Sym_n,$ symmetric matrices, 
\item[] $\Sk_n,$ $n=2k$, skew-symmetric matrices.
\end{itemize}
Our principal aim is to derive a classification of simple skew-symmetric families from the known simple classification of arbitrary square families. On some occasions we will make comparisons with the symmetric case too.  

\medskip
The simple families we are looking for are with respect to the following equivalences.
 
In the $\Sq_n$ case, we say that two matrix families, $M_1$ and $M_2,$ are  {\em $\Sq$-equivalent} if there exist a biholomorphism $\f$ of the source and two holomorphic maps $A,B: U \to \GL_n(\CC)$ such that
\begin{equation}\label{Esqeq}
M_1 \circ \f = A^\top M_2 B\,.
\end{equation}
For the $\Sym$- and {\em $\Sk$-equivalences},
we require the existence of only one holomorphic $s$-parameter family $A$ of invertible $n \times n$ matrices such that
\begin{equation}\label{Eskeq}
M_1 \circ \f = A^\top M_2 A\,.
\end{equation}
We use a subscripted version of the terminology, like $\Mat_n$-equivalence, if we want to emphasise the size of the matrices. 

Germ versions of the above definitions are straightforward. 
An equivalence class of germs of matrix families will be called a {\em matrix singularity}.  We denote by $\O_s$ the space of all holomorphic function germs on $(\CC^s,0),$ and identify the space of all holomorphic map germs $M:  (\CC^s,0) \to \Mat_n$ with the $\O_s$-module $\O_s^N$ where $N= {\rm dim}_\CC\Mat_n.$  We write $\O_x$ instead of $\O_s$ when we want to emphasise the coordinates $x$ introduced on $\CC^s.$ As usual, $\frak m_x$ will be the maximal ideal in $\O_x.$

Our matrix classification will make use of right ($\R$) and contact ($\K$) equivalence of holomorphic functions, and of their versions $\R_\p$ and $\K_\p$
for functions on manifolds with boundary \cite{Ab}.
In the $\Sq_1$, $\Sym_1$ and $\Sk_2$ cases, our equivalences are the same as the $\K$-equivalence of functions, and, therefore, our classification theorems do not mention these lowest size cases. The $\K$-equivalence class of the determinantal or respectively Pfaffian function of a matrix family is an invariant of the $\Mat$-equivalence.

\medskip
Let $E_{j\ell}$ be the $n\times n$ matrix with the $j\ell$-entry 1, and all other entries zero. The extended tangent spaces to the 
equivalence classes of germs  $M:  (\CC^s,0) \to \Mat_n,$ with the source coordinates $x_1,\dots,x_s,$ are 
$$
\begin{array}{rcl}
T_{\Sq}M & = & \O_s \left<\p M/\p x_i, \   E_{j\ell}M, \  ME_{pq} \right>_ { i=1,\dots,s; \   j,\ell,p,q=1,\dots,n},  \qquad {\rm and}
\\ \vspace{-5pt} \\
T_{\Mat}M & = & \O_s \left<\p M/\p x_i, \   E_{j\ell}M + ME_{\ell j} \right>_{ i=1,\dots,s; \   j,\ell=1,\dots,n}
\quad
{\rm for} \ \ \Mat=\Sym,\Sk.
\end{array}
$$
We will denote by $\tau_{\Mat}(M)$ 
the {\em Tjurina number\/} of $M,$ that is, the codimensions of the above extended tangent spaces in $\O_s^N.$ When dealing with skew-symmetric matrices we will write $\E_{ij}$ for $E_{ij}-E_{ji}.$

\medskip 
All our groups of matrix equivalences are in Damon's class of geometric subgroups of the $\K$-equi\-valence group (see \cite{D}), and hence a $\Mat$-miniversal deformation of a matrix family germ $M$ may be written as
$$
M + \l_1 \f_1 + \dots + \l_\tau \f_\tau\,,
$$ 
where $\tau = \tau_{\Mat} (M),$ 
and the $\f_i \in \O_s^N$ form a basis of the quotient $\O_s^N/T_{\Mat} M.$ 

\medskip
We call a matrix family  $M:  (\CC^s,0) \to (\Mat_n,0)$ {\em quasi-homogeneous\/} if it is possible to assign positive weights $w_1, \dots, w_s$ to the source coordinates so that each entry $m_{ij}(x_1,\dots,x_s)$ would be a quasi-homogeneous function of degree $d_{ij}$ and each $2 \times 2$ minor would also be quasi-homogeneous, that is, $d_{ij} + d_{pq} = d_{iq} + d_{pj}$ must hold for the entire range of the indices. The last condition implies in particular that $\det\circ M$ or $\Pf\circ M$ are quasi-homogeneous too. All normal forms of matrix families appearing in this paper are quasi-homogeneous.

\bigskip 
By the {\em matrix corank\/} of a germ $M:  (\CC^s,0) \to \Mat_n$ we will understand the corank of the matrix $M(0).$ In our symmetric and arbitrary square settings, a family of the matrix corank $c$ is equivalent to a family 
$\left( \begin{array}{cc}
M' & 0 \\ 0 & I_{n-c}
\end{array}
\right)$ where $M'$ is a germ of a $c \times c$ matrix family of the same type and $M'(0)$ is the zero matrix.
Similar reduction exists in the skew-symmetric case, with the only difference that the identity corner should be replaced by the block-diagonal matrix $J_{n-c}$ with the elementary  blocks
$J_2=\left( \begin{array}{cc} 0 & 1 \\ -1 & 0 \end{array} \right)$ along the diagonal.

Two germs of matrix families $M_i: (\CC^{s},0) \to \Mat_{n_i},$ $i=1,2,$  will be called {\em stably $\Mat$-equivalent} 
if there exists $n$ such that the two `extended'  families
$\left( \begin{array}{cc}
M_i & 0 \\ 0 & I_{n-n_i}
\end{array}
\right)$, or respectively $                        
\left( \begin{array}{cc}
M_i & 0 \\ 0 & J_{n-n_i}
\end{array}
\right)$, are $\Mat_n$-equivalent. All our $\Mat$-simple classifications are up to the stable $\Mat$-equivalences. All classification tables will be for $M(0)=0.$


For a classification process, it is useful to have an a priori weaker intermediate notion of simplicity: 

\begin{definition} \label{Dmatc} {\em
We will call a map germ $M: (\CC^{s},0) \to (\Mat_c,0)$ (that is, that of the matrix corank $c$) {\em $\Mat_c$-simple} if it is simple within the set of all germs of the same matrix corank $c.$}
\end{definition}

 Adjacencies of such an $M$ to germs of lower matrix coranks are ignored by the $\Mat_c$-simplicity. This notion allows us to first find lists of $\Mat_c$-simple singularities for various  $s$ and $c,$ and only after that to check for the existence or non-existence of adjacencies to non-simple singularities of lower matrix coranks.  One of the results of the classifications carried out in \cite{B,BT} is

\begin{prop} \label{Pmatc}
In the symmetric and arbitrary square settings, a matrix family germ is $\Mat$-simple if and only if it is $\Mat_c$-simple, where $c$ is the matrix corank of the germ.
\end{prop}

Thus the lists of $\Mat$-simple singularities in each of these two cases is the union of all the $\Mat_c$-simple lists. We will show that the same is true in the skew-symmetric case as well.  

\

\section{Simple singularities of arbitrary square matrices} \label{Ssq}  
Our main Theorem \ref{Tmain} that classifies $\Sk$-simple matrix families of even size will be derived from the $\Sq$-simple classification obtained in \cite{BT}. We shall now recall that classification. The tables below contain normal forms of the singularities and their Tjurina numbers $\tau_{\Sq}.$ Whenever the determinantal functions of the families have isolated singularities, we point out the classes of these functions. We are saving the notation $s$ and $M$ for skew-symmetric matrix families, and therefore using $\cs$ and $\cM$ in the $\Sq$ case.

\begin{theorem} \label{TBT}  {\em \cite{BT}}
All  $\Sq$-simple germs $(\CC^\cs,0) \to (\Sq_k,0)$ 
are those that appear in the following list.

{\em (1)} When $\cs=1$ all finitely $\Sq$-determined germs are simple and $\Sq$-equivalent to a germ of the form $\diag(x^{a_1},x^{a_2}, \dots, x^{a_k})$ where $a_1 \le a_2 \le \dots \le a_k<\infty.$ This germ has $\tau_{\Sq} = \sum_{i=1}^k (2(k-i) +1)a_i - 1.$

{\em (2)} When the target corank of the derivative of the germ is 0 we have a normal form $\cM: \CC^{k^2}_x \times \CC^{\cs-k^2}_z \to \Sq_k$ given by $m_{ij} = x_{ij}$ which is $\Sq$-simple, and has $\tau_{\Sq} =0.$ (The $z$s are redundant variables.)

{\em (3)} When the target corank of the derivative of the germ is 1 
we have two cases.

\ \ \ \ \   {\em (i)} A normal form $\cM: \CC^{k^2-1}_x \times \CC^{\cs-k^2+1}_z \to \Sq_k$ given by 
$$
\left( \sum\limits_{i=2}^k x_{ii} + f(z) \right) E_{11} + \sum \limits_{(ij)\ne(11)} x_{ij} E_{ij}
$$
where $f$ is one of Arnold's $\R$-simple germs $A_\mu, D_\mu, E_6, E_7, E_8.$ In this case $\tau_{\Sq}(\cM)$ coincides with the Tjurina number number of $f.$

\ \ \ \ \   {\em (ii)} A normal form $\cM: \CC^{k^2-1}_x \times \CC^{\cs-k^2+1}_z \to \Sq_k$ given by 
$$
\left( \sum\limits_{i=2}^{k-1} x_{ii} + f(x_{kk},z) \right) E_{11} + \sum \limits_{(ij)\ne(11)} x_{ij} E_{ij}
$$
where $f$ is one of Arnold's $\R_\p$-simple germs $B_\mu, C_\mu, F_4$ on $\CC^{\cs-k^2+2}_{x_{kk},z}$ with the boundary $x_{kk}=0$.
Here $\tau_{\Sq}(\cM)$ coincides with the $\R_\p$ Tjurina number of $f.$

{\em (4)} {\em (see also \cite{Gvc})} If the target corank of the derivative of the germ is $2$ while $k >2$ then the only possibility is $k=3$ with $\cs=7,$ in which case the $\Sq$-simple classification consists of the following 
series and three sporadic singularities, whose $\tau_{\Sq}$ are respectively  $p+1, 4, 5$ and $6:$   
$$
\begin{array}{cc}
\left( 
\begin{array}{ccc}
-x_{22} & x_{12} & x_{13} \\
- x_{33} + x_{12}^p & x_{22} & x_{23} \\
x_{31} & x_{32} & x_{33}
\end{array}
\right)\!\!, p \ge 1, 
 &
\left( \begin{array}{ccc}
-x_{22}-x_{33} & x_{12} & x_{13} \\
-x_{32}+x_{13}^2 & x_{22} & x_{23} \\
x_{31} & x_{32} & x_{33}
\end{array}
\right) , 
\\
&
\\
\left( \begin{array}{ccc}
-x_{22}-x_{33} + x_{13}^2 & x_{12} & x_{13} \\
-x_{32} & x_{22} & x_{23} \\
x_{31} & x_{32} & x_{33}
\end{array}
\right) ,  &  
\left( \begin{array}{ccc}
-x_{22}-x_{33} & x_{12} & x_{13} \\
-x_{32}+x_{13}^3 & x_{22} & x_{23} \\
x_{31} & x_{32} & x_{33}
\end{array}
\right) . 
\end{array}
$$

{\em (5)} When $k=2$ the simple germs that are not covered in {\em (1)--(3)} are given in Table $1$  and Table $2.$ Matrices $
\left( \begin{array} {cc}  a & b \\ c & d \end{array} \right)$ are written there in the line format $(a,b,c,d).$

\begin{table} [h] \label{T2to2} 
\caption{$\CC^2 \to \Sq_2.$}
\vspace{-16pt}
$$\begin{array}{l|c|c||l|c|c}
\hline
\vspace{-10pt} &&&&& \\
Normal \  form & \det & \tau_{\Sq} & Normal \  form & \det & \tau_{\Sq}
\\
\vspace{-12pt} &&&&& \\
\hline \hline
\vspace{-7pt} &&&&& \\
(x,y^p,y^q,x), \,  1\le p \le q & A_{p+q-1} & 2p\!+\!q\!-\!1 &
(x,0,0,xy+y^p), 3\le p & D_{2p} & 3p \\
\!\!\!\left. \begin{array}{l}(x,y,x^2+y^p,0) \\  (x,x^2+y^p,y,0) \end{array}
\right\}\!\!, \,   2\le p & D_{p+2} & p+3 &
(x,y^p,y^q,xy), \  \underbrace{2 \le p,q} & D_{p+q+1}
& \!\! \begin{array} {l} p\!+\!q\!+\!1 
\vspace{-4pt}
\\ +\! \min\{p,q\}\end{array}\!\!\!\!\!
\\
\!\!\!\left. \begin{array}{l}(x,y,y^3,x^2) \\  (x,y^3,y,x^2) \end{array}
\right\}  & E_6 & 7 &
(x,y^2,y^2,x^2) & E_6 & 8 \\
\!\!\!\left. \begin{array}{l}(x,y,xy^2,x^2) \\  (x,xy^2,y,x^2) \end{array}
\right\}  & E_7 & 8 &
\!\!\!\left. \begin{array}{l}(x,y^2,0,x^2+y^3) \\  (x,0,y^2,x^2+y^3) \end{array}
\right\}  & E_7 & 9 \\
\!\!\!\left. \begin{array}{l}(x,y,y^4,x^2) \\  (x,y^4,y,x^2) \end{array}
\right\}  & E_8 & 9 &
(x,0,0,x^2+y^3) & E_7 & 10 \\
(x,0,0,y^2+x^p), 2 \le p & D_{p+2} & p+4 &
\!\!\!\left. \begin{array}{l}(x,y^2,y^3,x^2) \\  (x,y^3,y^2,x^2) \end{array}
\right\}  & E_8 & 10 \vspace{-12pt} \\
&&&&& \\
\hline
\end{array}$$
\end{table}

\begin{table} [h] \label{T3to2} 
\caption{$\CC^3 \to \Sq_2.$}
\vspace{-4pt}
$$\begin{array}{l|c||l|c}
\hline
\vspace{-10pt} &&& \\
Normal \  form & \det & 
Normal \  form & \det 
\\
\vspace{-12pt} &&& \\
\hline \hline
\vspace{-7pt} &&& \\
(x,z^p,z^q,y), \  1\le p \le q & A_{p+q-1} & 
(x,y,y^2+z^3,x) & E_7 
\\
(x,y,z^2+y^p,x), \  2 \le p & D_{p+2} & 
(x,y,yz,x+z^p), \  2 \le p & D_{2p+1} 
\\
(x,y,y^2,x+z^2) & E_6 & 
(x,y,yz+z^p,x), \  3 \le p & D_{2p} 
\vspace{-12pt} \\
&&& \\
\hline
\end{array}$$
\end{table}

{\em (6)} When $k=3$ the simple germs that are not covered in {\em (1)--(4)} are given in Table $3.$ 

\begin{table} [h] \label{T2to3}  
\caption{$\CC^2 \to \Sq_3.$}
\vspace{-4pt}
$$\begin{array}{l|c|c||l|c|c}
\hline
\vspace{-10pt} &&&&& \\
Normal \  form & \det & \tau_{\Sq}, \  \tau_{\Sk}^{triv} & Normal \  form & \det & \tau_{\Sq}, \  \tau_{\Sk}^{triv}
\\
\vspace{-12pt} &&&&& \\
\hline \hline
\vspace{-7pt} &&&&& \\
\begin{array}{c}
\left(\begin{array}{ccc}
x & y^p & 0 \\
y^q & x & 0 \\
0 & 0 & y
\end{array}
\right) \\  1\le p \le q 
\end{array} & D_{p+q+2} & 
\begin{array}{l}
2p+q+4, \\
 3p+q+8 
\end{array} &
\left(\begin{array}{ccc}
x & y & 0 \\
0 & x & y \\
y^2 & 0 & x
\end{array}\right) & E_6 & 9, \  15
\\ &&&&& \\
\left(\begin{array}{ccc}
x & y & 0 \\
0 & x & y \\
xy & 0 & x
\end{array}\right) & E_7 & 10, \  16  &
\left(\begin{array}{ccc}
x & y & 0 \\
y^2 & x & 0 \\
0 & 0 & x
\end{array}\right) & E_7  & 11, \   19
\\ &&&&& \\
\left(\begin{array}{ccc}
x & y & 0 \\
0 & x & y \\
y^3 & 0 & x
\end{array}\right) 
& E_8  &   11, \   17     &
\left(\begin{array}{ccc}
x & y & 0 \\
0 & x & y^2 \\
y^2 & 0 & x
\end{array}\right) & E_8  & 12, \  20
\vspace{-12pt} \\
&&&&& \\
\hline
\end{array}$$
\end{table}

\end{theorem}

\begin{rems} \label{RBT}  {\em
(i) Members of each braced up pair of mutually transposed normal forms in Table 1 are not $\Sq$-equivalent.
The $p \leftrightarrow q$ pairing within the $D_{p+q+1}$ series is of the same nature.
We point out in the next section the reason why such pairs give rise to $\Sk$-equivalent singularities. 

(ii) Non-symmetric normal forms of all other singularities are $\Sq$-equivalent to their transposes. 

(iii) For  all maps of Table $2,$ their Tjurina numbers coincide with the Milnor numbers of their determinantal functions (see \cite{GM} for a general case).

(iv) There is some repetition in the theorem. Namely, the $1=p\le q$ subseries of the first series in Table 2 is already covered by part (3). This duplication is permitted for the convenience of later references.

(v) The meaning of the Tjurina numbers $\tau_{\Sk}^{triv}$ in Table 3 will be clarified by Theorem \ref{Tmain}.}
\end{rems}

\ 

\section{Correspondence between $\Sq_k$- and $\Sk_{2k}$-equivalences} \label{Scor}

\subsection{Suspensions} \label{Ssus}
The classification of simple $2k \times 2k$ skew-symmetric families $M$ turns out to be directly related to that of $k \times k$ arbitrary square families $\check M$ mostly by two kinds of suspensions which split $M$ into four $k \times k$ blocks:

\begin{definition} {\em 
The skew-symmetric matrix family $M = \left( \begin{array} {cc} 0 & \check M \\ - \check M^\top & 0 \end{array} \right)$ will be called the {\em trivial suspension of \/} a square matrix family $\check M.$ The {\em full suspension of \/} $\check M$ is the family $M = \left( \begin{array} {cc} V & \check M \\ - \check M^\top & W \end{array} \right),$ where the blocks $V$ and $W$ are given by germs of diffeomorphism of two different copies of $(\CC^{k(k-1)/2},0)$ with $(\Sk_k,0).$}
\end{definition}

Thus, the number of parameters in the trivial case stays the same, and in the full case it increases by $k(k-1).$ 
We take the set of all upper-triangle entries of the blocks $V$ and $W$ for the additional parameters of a fully suspended map.

\medskip
Let $\D_k$ be the subgroup of the $\Sk_{2k}$-equivalence group which preserves the 4-block structures of the trivial suspensions. The right group of $\D_k$ is that of diffeomorphisms of the domain --- let it be $(\CC^\cs,0)$ --- of $\check M.$  The group of row-column operations 
allowed in $\D_k$ has two connected components. In terms of (\ref{Eskeq}), its identity component consists of all
matrix families 
$
A = \left( \begin{array} {cc} A_1 & 0 \\ 0 & A_2 \end{array} \right)
$,
$A_1,A_2: (\CC^\cs,0) \to \GL_k(\CC)$,
acting on the trivial suspensions by
$$
A^\top \left( \begin{array} {cc} 0 & \check M \\ - \check M^\top & 0 \end{array} \right)  A = 
\left( \begin{array} {cc}
0  & A_1^\top \check M A_2 \\
- (A_1^\top \check M A_2)^\top & 0
\end{array}
\right).
$$
The other connected component contains the constant family $A =  \left( \begin{array} {cc} 0 & I \\ -I & 0 \end{array} \right)$ that sends $\check M$ to  $\check M^\top.$  

\begin{definition} \label{Dext} {\em 
The $\Sq$-equivalence of arbitrary square matrix families which additionally considers a pair of mutually transposed families as equivalent will be called the} 
$Sq^{\!\top}$-equivalence.
\end{definition}

This is exactly the equivalence of  the families $\check M$ up to the action of the group $\D_k.$
Thus we have:

\begin{prop} \label{Ptriv}
The $\Sk$-classification of those maps $M: (\CC^\cs,0) \to (\Sk_{2k},0)$ which are $\Sk$-equivalent to trivial suspensions of maps $\check M: (\CC^\cs,0) \to (\Sq_k,0)$ coincides with the $\Sq^{\!\top}$-classification of the maps $\check M$.
\end{prop}

A similar observation turns out to be valid for the full suspensions as well. To see this, we consider the same group $\D_k,$ whose right group is still that of diffeomorphisms of $(\CC^\cs,0).$ The action of the $\GL_k \times \GL_k$ part of $\D_k$ will now be compensated by changes of the suspending coordinates. Such changes may depend on the coordinates on $\CC^\cs.$

Namely, the results of the actions
of $
\left( \begin{array} {cc} A_1 & 0 \\ 0 & A_2 \end{array} \right)\!,\, \left( \begin{array} {cc} 0 & I \\ -I & 0 \end{array} \right) \in \D_k$ on a full suspension of $\check M$ are respectively
\begin{equation} \label{ED}
\left( \begin{array} {cc}
A_1^\top V A_1 & A_1^\top \check M A_2 \\
- (A_1^\top \check M A_2)^\top & A_2^\top W A_2
\end{array}
\right) \quad {\rm and}  \quad 
\left( \begin{array} {cc} W & \check M^\top \\ - \check M & V \end{array} \right).
\end{equation}
In the first case here, the entries of the blocks $A_1^\top V A_1$ and $A_2^\top W A_2$ 
provide replacements for the entries of $V$ and $W$  as coordinates on the suspended source. In the second case, we can just swap $V$ and $W$. This gives us

\begin{prop} \label{Pfull}
The $\Sk$-classification of those maps $M: (\CC^{\cs+k(k-1)},0) \to (\Sk_{2k},0)$ which are $\Sk$-equivalent to full suspensions of maps $\check M: (\CC^\cs,0) \to (\Sq_k,0)$ coincides with the $\Sq^{\!\top}$-classification of the maps $\check M$.
\end{prop}

Thus, our full suspensions are similar to function suspensions (or stabilisations) $z^2 + f(x)$ by the square of a new variable provided by the Morse Lemma with parameters which uses fibred coordinate changes $(Z(z,x),X(x)).$

\bigskip
In our classification, we will also need a third suspension type, in the $\Sk_4$ situation only:

\begin{definition} \label{Dint} {\em
The {\em intermediate suspension} of a map $\check M: (\CC^\cs_x, 0) \to (\Sq_2,0)$ is the map $M: (\CC^{\cs+1}_{x,z}, 0) \to (\Sk_4,0)$: \quad
$M = \left( \begin{array} {cc} 0 & \check M \\ - \check M^\top & 0 \end{array} \right)+ zJ_4.$}
\end{definition}

\begin{rems} \label{RPf} {\em (i) From (\ref{ED}), we notice that the group $\D_2'$ preserving this suspension shape is a subgroup of $\D_2$ singled out by the coincidence of the determinantal functions of the families $A_1$ and $A_2.$ 
However, the intermediate suspension will appear in what follows just in the quasi-homogeneous context, and there is no difference between the $\D_2$- and $\D_2'$-classifications of the maps $\check M$ which have quasi-homogeneous representatives (see \cite{Gvc}).

(ii) For the intermediate suspensions,  $\Pf \circ M = z^2 - \det \circ \check M.$ For the full suspensions of the $\Sq_2$-families, $\Pf \circ M = vw - \det \circ \check M.$

(iii) Trivial suspensions in the $k=2$ case were considered in \cite{H}.}
\end{rems}

\subsection{$\Sk$-simple classification} \label{Smain}
The main result of this paper is that the list of all $\Sk$-simple matrix families may be produced by appropriate suspensions from the list of all $\Sq$-simple singularities of Theorem \ref{TBT}:

\begin{theorem} \label{Tmain}
There are bijections between the sets of $\Sk$-simple classes  of map germs $M: (\CC^s,0) \to (\Sk_{2k},0)$ 
and the sets of $\Sq^{\!\top}$-simple classes of map germs $\check M: (\CC^{\check s},0) \to (\Sq_{k},0).$ The bijections are established by trivial, full and intermediate suspensions. The types of the suspension correspondence depend on $s$ and $k$ as shown in Table $4$.
\begin{table}[h]
\caption{Suspension correspondence between $\Sk$- and $\Sq$-simple matrix singularities}
$$\begin{array}{c|c|c|c|c}
\hline
\vspace{-10pt} &&&& \\
k  &  s  & \check s & suspension &
\begin{array}{c} \vspace{-4pt}
related \  part \\
of \  Theorem \  \ref{TBT}
\end{array}
\\
\vspace{-12pt} &&&& \\
\hline \hline
\vspace{-7pt} &&&& \\
any  & 1 & 1 & trivial &  \rm{(1)} \\
1 & any & s & trivial \equiv full \\ 
2 & 2 & 2 & trivial  & \rm{(5), \  Table \  1} \\
2 & 3 & 2 & intermediate & \rm{(5), \  Table \  1} \\
2 & \ge 4 & s-2 & full & \rm{(5,2,3)} \\ 
3 & 2 & 2 & trivial & \rm{(6)} \\
3 & 13 & 7 & full & \rm{(4)} \\
\ge 3 & \ge k(2k-1) - 1 & s-k(k-1) & full & \rm{(2,3)}
\vspace{-12pt} \\
&&& \\
\hline
\end{array}
$$
\end{table}
\end{theorem}

By Remarks \ref{RBT}(i,ii),  the only difference between the $\Sq$- and $\Sq^\top$-simple lists is that singled out by the braces in Table 1. 

\medskip
Suspension-wise, relations between the Tjurina numbers $\tau_{\Sk} (M)$ and $\tau_{ \Sq}(\cM)$, and between miniversal deformations $\M$ and $\check \M$ of the corresponding simple families are as follows:

\smallskip
({\em full}) In all such cases, $\tau_{\Sk} (M) = \tau_{ \Sq}(\cM).$ One can take for $\M$ the simultaneous full suspension $\left( \begin{array} {cc} V & \check \M \\ - \check \M^\top & W \end{array} \right)$ of $\check \M.$

\smallskip
({\em intermediate})  This type appears in Table 4 only once, when $k=2, s=3$ and $\cs = 2$. We then have  $\tau_{\Sk} (M) = \tau_{ \Sq}(\cM) + \kappa,$ where $\kappa$ is the codimension in $\O_{2}$ of the ideal $I_{\cM}$ generated by all four entries of the $2 \times 2$ matrix $\cM.$ 
It is possible to set
$$\M = \left( \begin{array} {cc} 0 & \check \M \\ - \check \M^\top & 0 \end{array} \right)+ zJ_4 + (\nu_1 \psi_1 + \dots + \nu_\kappa \psi_\kappa) \E_{12},$$
where the $\nu_i$ are the additional deformation parameters and the $\psi_i$ represent a basis of $\O_{2}/I_{\cM}.$

\smallskip
({\em trivial}) The best overall description we have is for the case $k=2, s=\cs = 2.$ Keeping the notation used just above, we now have $\tau_{\Sk} (M) = \tau_{ \Sq}(\cM) + 2\kappa$ and can take
$$\M = \left( \begin{array} {cc} 0 & \check \M \\ - \check \M^\top & 0 \end{array} \right)+ (\nu_1 \psi_1 + \dots + \nu_\kappa \psi_\kappa) \E_{12} +   (\o_1 \psi_1 + \dots  + \o_\kappa \psi_\kappa) \E_{34},$$
where the $\o_i$ provide a further extension of the set of deformation parameters.

For $s= \cs =1$ and $k$ arbitrary,  $\tau_{\Sk} (M) = 2\tau_{ \Sq}(\cM) + 1.$ A way to construct $\M$ from $\check \M$ will be clear from Section  \ref{Sfrom1}.

Finally, the behaviour of the $k=3,$ $s = \cs =2$ singularities is idiosyncratic and we gave their Tjurina numbers in Table 3.

\medskip
All the relations between the Tjurina numbers listed here will be explained in Section \ref{Sred} during detailed discussions of the particular situations. 

\begin{rem} {\em 
All normal forms of matrix families in Theorem \ref{Tmain} are quasi-homogeneous. Therefore, the list of simple matrix singularities stays the same for a finer version of the $\Sk$-equivalence when the maps $A$ in  (\ref{Eskeq}) take values in $\SL_n(\CC)$ (cf. \cite{GM, Gvc}).}
\end{rem}

\

\section{Dimensional constraints for $\Sk$-simplicity} \label{Sdim}
This section starts our proof of Theorem \ref{Tmain}. 
We derive constraints on the linear parts of $\Sk$-simple matrix families $(\CC^s,0) \to (\Sk_{2k},0)$ imposed by the restriction of the $\Sk$-equivalence (\ref{Eskeq}) to the 1-jets of the families. Therefore, we will be mainly focused on the action
\begin{equation}\label{Elinact}
(A,X) \mapsto A^\top X A, \quad A \in \GL_{2k}, \ X \in \Sk_{2k},
\end{equation}
and the induced actions of $\PGL_{2k}$ on the Grassmannians $\Gr^r_N$ of $r$-dimensional subspaces  
of $\Sk_{2k}.$ All through this Section $N=k(2k-1).$

One of our tools will be 
a skew-symmetric version of the dualities used in \cite{B, BT}. Namely, the space $\Sk_{2k}$ is self-dual 
with respect to the non-degenerate quadratic form $\left< X, Y \right> = \tr(XY).$ 
For an element $Z \in \Gr^r_N,$ there is the element $Z^\perp = \{Y \in \Sk_{2k} | \tr(XY)=0, \  \forall X \in Z\} \in \Gr^{N-r}_N.$ We have

\begin{prop} \label{Pdual}  {\em (cf. \cite{B,BT}).}
Let $Z_i,$ $i \in I,$ be a listing of the subspaces of $\Sk_{2k}$ of dimension $r$ up to the $\GL_{2k}$-equivalence. Then the subspaces $Z_i^\perp,$ $i \in I,$ is a listing of the subspaces of dimension $N-r$
in $\Sk_{2k}$ up to the same equivalence. 
\end{prop}

Indeed, orbit-wise we have for $A \in \GL_{2k}$ and $Z \in \Gr^r_N$: \quad 
$
(A^\top Z A)^\perp = A^{-1} Z^\perp (A^{-1})^\top.
$

\bigskip 
Let $r$ now be the rank of a map germ $M:(\CC^s,0) \to (Sk_{2k},0),$ $s \ge r,$ at the origin. What conditions must $r$ and $k$ meet for $M$ to have a chance to be $\Sk$-simple? 

The image of the derivative of $M$ at the origin is an element $\a \in \Gr_N^r.$ Therefore, one of necessary conditions for the simplicity is the openness of the generic $\PGL_{2k}$-orbit in $\Gr^r_N.$ For this, the dimension of the Grassmannian must at least not exceed the dimension of the group, that is,
\begin{equation}\label{Egr}
r (N- r) \le  4k^2  - 1. 
\end{equation}
This inequality holds for $r=0,1,2,N-2,N-1,N$ with any $k,$ and allows some extras for two exceptional sizes: 
\begin{itemize}
\item[] $k=2$ with $r=3 = N-3,$ and
\item[] $k=3$ with $r=3, 12$ (here $3+12 = N$).
\end{itemize}

To reduce the range obtained, 
consider the restriction of the Pfaffian function from $\Sk_{2k}$ to $\a.$ This restriction is a $k$-form  on $\CC^r.$ Its contact class is invariant under the equivalences (\ref{Eskeq}) and (\ref{Elinact}). Therefore, it has moduli if $k \ge 4$ and $r=2,$ or $k=3$ and $r=3.$ (These moduli are easily seen to be realisable by variations of $\a.$) This also  means that generic orbits of the $\PGL_{2k\ge 8}$-actions on the $\Gr^2_N$ and of the $\PGL_{6}$-action on $\Gr^3_6$ are not open. Proposition \ref{Pdual} implies similar non-openness for the $\Gr^{N-2}_N$ and $\Gr^{12}_6$. So we have:

\begin{prop}\label{Prank}
Simple skew-symmetric matrix families may exists only in the following cases: 
\begin{itemize}
\item[] $r=0,1,N-1,N$ for any $k$;
\item[] $k=2$ for $r=2,3,4$;
\item[] $k=3$ for $r=2,13.$ 
\end{itemize}
\end{prop} 

In particular, in $\Sk_4,$ this allows any $r \le 6 = N$.

\medskip
If the dimension $s$ of the domain of a map $M$ is greater than $r$ then $M$  deforms to maps of higher ranks at the origin.
Thus, Proposition \ref{Prank} implies that there are no $\Sk$-simple families in $2\le s \le N-2$ parameters if $k\ge 4,$
and in $3 \le s \le 12$ parameters if $k=3.$ This yields

\begin{corollary} \label{Cdim}
$\Sk$-simple matrix families may exist only for the following dimensional triplets $s, k, r$:
\begin{itemize}
\item[a)] $s=1,$ $r=0,1$ and any $k;$ 
\item[b)] $s \ge N-1,$ $r=N,N-1$ and any $k\ge 3;$
\item[c)] $k=3,$ $s=2 \ge r$;
\item[d)] $k=3,$ $r=13 \le s$;
\item[e)] $k=2,$ $s \ge 2.$
\end{itemize}
\end{corollary}

\section{Reduction to the suspensions} \label{Sred}
We shall now analyse the sequence of dimensional options allowed by Corollary \ref{Cdim}. Until Section \ref{Sskcheck},  our aim will be 
to detect $\Sk_{2k}$-simple families (see Definition \ref{Dmatc}), that is, 
we shall not be interested in adjacencies between singularities of different matrix coranks.  In each case $\Sk_{2k}$-simple singularities will be reduced to the suspension form promised by Theorem \ref{Tmain}.
The $\Sk_2$ case is trivial and therefore will be omitted.

Finally, in Section \ref{Sskcheck}, we bring all our fixed-matrix-corank classifications together and prove that none of our $\Sk_{2k}$-simple singularities is adjacent to non-simple singularities of lower matrix corank.

We denote by $T_{ij}$ the operation of addition of the $i$th column of a matrix to its $j$th column accompanied by the same transformation of the rows. We write $\a T_{ij}$ if the $i$th column and row are used in this context with a coefficient $\a.$ For example, $xT_{11}$ multiplies the first row and the first column by $1+x.$ 

Block-diagonal matrices with an ordered sequence of square blocks $A_1, A_2, \dots, A_m$ of various size along the diagonal will be written as $A_1 \oplus A_2 \oplus \dots \oplus A_m.$

 \subsection{One-parameter families} \label{Sfrom1}

\begin{prop} \label{Ps1} {\em \cite{H}}
Any map germ $M: (\CC_x,0) \to (\Sk_{2k},0)$ with a finite $\Sk$ Tjurina number
reduces to the $2 \times 2$ block-diagonal form 
$$x^{a_1}J_2 \oplus x^{a_2}J_2 \oplus \dots \oplus x^{a_k}J_2, \quad
1 \le a_1 \le a_2 \le \dots \le a_k.$$ 
Such a map has $\tau_{\Sk} = \sum\limits_{i=1}^k 2a_i(2(k-i)+1) -1.$
\end{prop}

\medskip
{\em Proof.} 
Let $a_1$ be the lowest order in $x$ across all entries of our matrix $M.$
Up to a permutation of rows  (and the same permutation of columns) of  $M,$ we may assume that $a_1$ is the order of the entry $m_{12}.$  Dividing the first row and the first column of $M$ by a unit, we make $m_{12}= x^{a_1}.$ A sequence of multiples of the operations $T_{2i},$ $i>2,$ followed by multiples of the $T_{1j},$ $j>2,$ reduces to zero all the entries in the first two rows and columns, except for $m_{12}$ and $m_{21}=-m_{12}.$  This brings $M$ to the block-diagonal form consisting of $x^{a_1} J_2$ and a $(2k-2) \times (2k-2)$ block $M_1.$ The claim of the proposition follows now by induction. 

The $\tau_{\Sk}$ formula is obvious.
\hfill{$\Box$} 

\medskip
All map germs with finite Tjurina numbers are, of course, simple in this case. The normal forms guaranteed by the proposition are clearly $\Sk_{2k}$-equivalent to the trivial suspensions of the $\Sq_k$-simple singularities $\cM=\diag(x^{a_1}, x^{a_2}, \dots, x^{a_k}).$

\subsection{Maps to $\Sk_4$} \label{Sto4}
In the subsections below, we are increasing the source dimension starting from $s=2.$
While the singularity of the Pfaffian of $M$ stays isolated, that is $s \le 6,$  our main $\Sk_4$-simplicity cutoffs will be by not allowing the Pfaffian to have singularities either of the {\em fencing} stable $\K$-equivalence classes 
$$
\begin{array} {rll}
P_8: & x^3 + y^3 + z^3 + axyz , \ \  &  a^3 + 27 \ne 0,\\
X_9: & x^4 + ax^2y^2 + y^4, &  a^2 \ne 4, \\
J_{10}: & x^3 + a xy^4 + y^6, &  4a^3 + 27 \ne 0,
\end{array}
$$
or adjacent to any of these classes. Fencing refers to the fact that any other non-simple singularity is adjacent to at least one element of this set, and this set is minimal for which such a property holds.

\subsubsection{Maps from $\CC^2$} \label{S2to4}
If a matrix family $M:(\CC^2,0) \to (\Sk_4,0)$  is $\Sk_4$-simple, then at least one of its entries must have non-trivial linear part. Otherwise $\Pf \circ M$ will in general be an $X_9$ function singularity and the $X_9$ modulus will be realisable in small perturbations of $M.$

Up to a permutation of rows and columns, we assume that $m_{14}$ is such an entry. So we take it for one of the coordinates on the source, and thus consider
$$
M = \left( \begin{array} {cccc}
0 & \a & \b & x \\
-\a & 0 & \g & \d \\
-\b & -\g & 0 & \e \\
-x & -\d & -\e & 0 
\end{array}
\right), \qquad {\rm where \  } \a,\b,\g,\d,\e \in \frak m_{x,y}.
$$             
A multiple of the $T_{42}$ move
allows us to reduce $\a$ to a function just in $y$. 
Similarly, multiples of $T_{43},T_{12}$ and $T_{13}$ eliminate  $x$ from $\b, \d$ and $\e$, so that we to arrive at
\begin{equation}
\label{E4x4}
 \left( \begin{array} {cccc}
0 & \a(y) & \b(y) & x \\
-\a(y) & 0 & \g(x,y) & \d(y) \\
-\b(y) & -\g(x,y) & 0 & \e(y) \\
-x & -\d(y) & -\e(y) & 0 
\end{array}
\right), 
\end{equation}
where modified entries are denoted as earlier. Up to a permutation of columns (and
in the same way of rows) and changes of signs of some entries, we can assume that $\b$ has the lowest order in $y$ comparing with $\a, \d$ and $\e.$ Then multiples of $T_{32}$ and $T_{14}$
make $\a =0$ and $\e=0.$ 
Thus, still keeping our earlier notation for modified entries, we end up with
\begin{equation} \label{E2to4}
M= \left( \begin{array} {cccc}
0 & 0 & \b(y) & x \\
0 & 0 & \g(x,y) & \d(y) \\
-\b(y) & -\g(x,y) & 0 & 0 \\
-x & -\d(y) & 0 & 0 
\end{array}
\right). 
\end{equation}

This is the trivial suspension of the matrix family $\cM: (\CC^2,0) \to (\Sq_2,0)$ written here in the top right corner.  All $\Sq_2$-simple matrix singularities in two variables (that is, those from Table 1) may be written in the form of this corner. The claim of Theorem \ref{Tmain} in relation to $\Sk_4$-simple maps from $\CC^2$ 
follows now from Propositions \ref{Ptriv} and \ref{Pmatc}. 

\medskip
The relation 
\begin{equation} \label{Etautriv}
\tau_{\Sk}(M) = \tau_{\Sq}(\cM) + 2\dim\, \O_2/I_{\cM}  
\end{equation}
from the end of Section \ref{Smain} comes from the observation that the moves $T_{32}, T_{42}, T_{31}, T_{41}$ add to the family (\ref{E2to4}) the matrix $\E_{12}= E_{12} - E_{21}$ multiplied by respectively $\beta, x, \gamma, \delta,$ and that we can similarly obtain the same multiples of $\E_{34}.$ Of course, the relation holds for an arbitrary dependence of $\cM$ on $x$ and $y,$ without  requiring one of the entries being exactly $x$ and two others depending on just one of the coordinates.

 \subsubsection{Higher-dimensional maps to $\Sk_4$} \label{Shto4}
A necessary condition for a germ $M: (\CC^s,0) \to (\Sk_4,0)$ to be either the intermediate or full suspension of a family of $2 \times 2$ matrices is that the rank of the quadratic part of the function $\Pf \circ M$ should be at least 1 or respectively at least 2 (see Remark \ref{RPf}ii). We will now show that the requirement of simplicity guarantees such ranks.

\begin{prop} \label{Pj1Pf}
Assume a family $M: (\CC^s,0) \to (\Sk_4,0),$ $s \ge 3,$ is $\Sk_4$-simple. Then {\em
\begin{itemize}
\item[(i)] {\em the rank of $j_0^2(\Pf \circ M)$ is positive; }
\item[(ii)] {\em moreover, this rank is at least 2 if $s \ge 4.$}
\end{itemize}}
\end{prop}

\medskip
{\em Proof.} (i) Up to a row-column permutation, the principal quasi-homogeneous part of a generic family $M$ failing property (i) is
$$
\left( \begin{array} {cccc}
0 & \ell_1 & \ell_2 & \ell_3 \\
-\ell_1 & 0 & q_3 & q_2 \\
-\ell_2 & -q_3 & 0 & q_1 \\
-\ell_3 & -q_2 & -q_1 & 0
\end{array} \right)
$$
where the $\ell_i$ and $q_j$ are respectively linear and quadratic forms on $\CC^s.$
The dimension of the space $\mathcal S$ of all such principal parts is $3(s + s(s+1)/2).$ The $\Sk_4$-equivalence group acts on $\SS$ by linear transformations of the source $\CC^s_{x_1,\dots,x_s},$ and by generating row-column operations $T_{11},$ $x_iT_{1j}$ and $T_{mj},$ where $i=1,\dots,s$ and $j,m=2,3,4.$ Taking the quasi-homogeneity into account, the dimension of the orbits of this action is at most
$s^2 + 3s + 9,$ which is smaller than $\dim \mathcal S$ if $s > 3.$

For $s=3,$ we can assume in general that all the $\ell_i$ are exactly the source coordinates $x_i.$ Then the 3-jet $x_1 q_1 - x_2 q_2 + x_3 q_3$ of the Pfaffian is in general a $P_8$ function singularity, and the modulus of the $P_8$ family is realisable by small changes in the $q_i.$
Thus, all matrix families $M$ failing property (i) are not $\Sk_4$-simple. 

\medskip
(ii) Property (i) implies that any of our $\Sk_4$-simple families may be written as an $(s-1)$-parameter deformation induced from an $\Sk$-miniversal deformation of the map $(\CC,0) \to (\Sk_4,0),$ $z \mapsto z J_4,$ that is,
$$
M = \left( \begin{array} {cccc}
0 & z + \l_1 & \l_2 & \l_3 \\
-z-\l_1 & 0 & \l_4 & \l_5 \\
-\l_2 & -\l_4 & 0 & z-\l_1 \\
-\l_3 & -\l_5 & -z+\l_1 & 0
\end{array} \right), \quad \l_1,\dots,\l_5 \in {\frak m}_{x_1,\dots,x_{s-1}}.
$$
Up to a row-column permutation, failure of part (ii) of  the proposition implies -- in the most generic situation -- a restriction on three of the terms here: $\l_1,\l_4,\l_5 \in {\frak m}^2_x.$ The principal quasi-homogeneous part of $M$ is singled out in this case by assigning the weights $w_z=3$ and $w_{x_i} =2$ to the variables. The quasi-degrees of the upper-triangle entries in the principal part are therefore $3, 2, 2; 4, 4; 3$ going along the rows and down. The space $\mathcal S'$ of all such principal parts has dimension $s^2 + s.$ Our $\Sk_4$-equivalence group acts on $\mathcal S'$ as $\GL_1 \times \GL_{s-1},$ and by the generating operations $T_{11},$ $x_i T_{12},$ $T_{22}$ and $T_{mj}$, where $i = 1, \dots , s-1$ and $m,j=3,4.$ The quasi-homogeneity implies that the orbits of this action have dimension at most $s^2-s+6,$ which is smaller than $\dim \mathcal S'$ when $s>3.$
\hfill{$\Box$}

\bigskip
Proposition \ref{Pj1Pf} reduces the $\Sk_4$-simple classifications to the $\Sq_2^\top$-simple classifications: 

\begin{corollary} \label{Chto4}
{\em (i)} For $s\ge 4,$ every $\Sk_4$-simple map germ $M:(\CC^s,0) \to (\Sk_4,0)$ is $\Sk$-equivalent to the full suspension of an $\Sq$-simple map
germ $\check M: (\CC^{s-2},0) \to (\Sq_2,0).$ This relation is a bijection between the $\Sk_4$-simple and $\Sq_2^\top$-simple classes of matrix families. 

{\em (ii)} Every $\Sk_4$-simple map germ $M:(\CC^3,0) \to (\Sk_4,0)$ is $\Sk$-equivalent to the intermediate suspension of an $\Sq$-simple map
germ $\check M: (\CC^{2},0) \to (\Sq_2,0).$ This relation is a bijection between the 
$\Sk_4$-simple  and $\Sq_2^\top$-simple classes of matrix families. 
\end{corollary}

{\em Proof.} (i) Following Section \ref{S2to4} and Table 1 of Theorem \ref{TBT}, there is just one $\Sk$-equivalence class of map germs from $(\CC^2,0)$ to $(\Sk_4,0)$ whose Pfaffian has a Morse singularity. This class is represented by the trivial suspension of
$\left( \begin{array} {cc}
y & z \\
z & y 
\end{array}
\right)$. However, we prefer to represent this singularity by $y J_2 \oplus z J_2.$
Part (ii) of Proposition \ref{Pj1Pf} implies that any $\Sk_4$-simple map germ $M:(\CC^s,0) \to (\Sk_4,0)$ may be considered as an $(s-2)$-parameter deformation induced from an $\Sk$-miniversal deformation of this special map, that is,
\begin{equation}\label{Es31}
M=
\left( \begin{array} {cccc}
0 & y & \mu_1 & \mu_2 \\
-y & 0 & \mu_3 & \mu_4 \\
-\mu_1 & -\mu_3 & 0 & z \\
-\mu_2 & -\mu_4 & -z & 0 
\end{array}
\right), \quad \mu_1, \dots, \mu_4 \in {\frak m}_{x_1,\dots,x_{s-2}}.
\end{equation}
This matrix is the full suspension of its $\mu$-corner. Part (i) of the corollary follows now from Proposition \ref{Pfull}.

\ 

(ii) By part (i) of Proposition \ref{Pj1Pf}, we can assume that $M$ is a 2-parameter deformation of the map $z \mapsto zJ_4$:
\begin{equation}\label{Es3}
M(x,y,z) =
\left( \begin{array} {cccc}
0 & z+\l_0 & \l_1 & \l_2 \\
-z-\l_0 & 0 & \l_3 & \l_4 \\
-\l_1 & -\l_3 & 0 & z- \l_0\\
-\l_2 & -\l_4 & -z+\l_0 & 0 
\end{array}
\right), \quad \l_0, \dots, \l_4 \in \frak m_{x,y}.
\end{equation}

If $j_0^1 \l_0 \ne 0$ then we are back to the situation (\ref{Es31}) in which the $\mu$-corner reduces this time to $
\left( \begin{array} {cc}
x^p & 0 \\
0 & x^q
\end{array}
\right)$, $1 \le p \le q.$
The maps obtained may be rewritten as the intermediate suspensions of the
 $
\left( \begin{array} {cc}
x^p & y \\
y & x^q
\end{array}
\right).
$

\medskip
Assume now $j_0^1\l_0=0$ in (\ref{Es3}). The $\Sk_4$-simplicity is going to help us to get rid of $\l_0$ completely.

To achieve that, we partially repeat 
the sequence of pre-normalising moves used in Section \ref{S2to4}. 
Namely, to avoid the Pfaffians in (\ref{Es3}) of the class $X_9$ or adjacent to that, one of the $\l_{>0}(x,y)$ terms must have a non-trivial linear part. Without any loss of generality we may assume $\l_2$ being such, and moreover $\l_2=x.$ After that, moves repeating those used earlier eliminate  $x$ from $\l_1$ and $\l_4.$ The two-step operation $(-aT_{42})\circ (aT_{13})$ with 
$a = (\l_0(x,y) - \l_0(0,y))/x$
removes $x$ from $\l_0$ too (amending only $\l_3$).

We can get rid of the remaining $\l_0(y)$ if its order is not strictly smaller than the orders of both $\l_1(y)$ and $\l_4(y)$. Indeed, the operations $T_{14} \circ T_{32}$ and $T_{41} \circ T_{32}$ change $\l_0$ by respectively $\l_1$ and $-\l_4,$ hence -- under this condition on the orders -- 
at least one of the compositions $(bT_{14}) \circ (bT_{32})$ and $(cT_{41}) \circ (cT_{32})$
is able to eliminate $\l_0$ with appropriate choice of the coefficients $b,c \in \O_y.$

The only other effect of  $T_{14} \circ T_{32}$ and $T_{41} \circ T_{32}$ on the whole matrix is the addition of respectively $-2\l_0$ to $\l_4$ and $2\l_0$ to $\l_1.$ Hence one of these operations will lower the order of either $\l_4$ or $\l_1$ to the order of $\l_0$ if the order condition of the previous paragraph has not been satisfied.

With $\l_0$ gone, the mapping $M$ in (\ref{Es3}) becomes the intermediate suspension of 
$\check M = 
\left( \begin{array} {cc}
\l_1(y) & x \\
\l_3(x,y) & \l_4(y)
\end{array}
\right).
$

In Section \ref{Ssus} we introduced the group $\D_k$ which preserved the full suspension shape of $\Sk_{2k}$ matrix families. Let $\D'_2$ be the subgroup of $\D_2$ singled out by the requirement that
 the determinantal functions of the blocks $A_1$ and $A_2$ used in (\ref{ED}) must coincide. The action of $\D'_2$ on the maps $\check M: (\CC^2,0) \to (\Sq_2,0)$ is finer than the $\Sq_2^\top$-equivalence provided by $\D_2.$
However, the list of $\Sq_2^\top$-simple maps is quasi-homogeneous, and hence it coincides with its $\D'_2$ version.  Actually, the simple list stays unchanged even under the more restrictive condition that $A_1$ and $A_2$ should be maps to $\SL_2(\CC)$ rather than to $\GL_2(\CC)$ (see \cite{Gvc}).   \hfill{$\Box$}
 
\ 

As for the formulas
$$
\tau_{\Sk}(M) = \tau_{\Sq}(\cM) + \dim\, \O_2/I_{\cM} \qquad  {\rm  and} \qquad   \tau_{\Sk}(M) = \tau_{\Sq}(\cM)$$ from the end  of Section \ref{Smain} for the Tjurina numbers of respectively the intermediate and full suspensions $M$ of $\Sq_2$-simple families $\cM,$ the first is very similar to (\ref{Etautriv}): the absence of the coefficient 2 now is due to the presence of the generator $\p M/ \p z = J_4 = \E_{12} + \E_{34}$ of the tangent space to the $\Sk$-equivalence class. The second formula here follows from the presence in the tangent space of the elements $\E_{12}$ and $\E_{34}$ on their own as the derivatives of $M$ with respect to the two suspending variables.

\subsection{Maps to $\Sk_6$}
By Corollary  \ref{Cdim}  we have two situations to discuss here: when the ranks of our maps $M$ are either 2 or 13, with the source of respectively dimension either 2 or at least 13. It will turn out that dimensions higher than 13 have no $\Sk_6$-simple singularities of rank 13.

\subsubsection{Maps of $\CC^2$}
We first notice that the rank at the origin of a simple map $M: (\CC^2,0) \to (\Sk_6,0)$
must be 2: had it been 1, then in the most general situation the function $\Pf \circ M$ would have had a $J_{10}$ singularity with the $J_{10}$ modulus realisable by small perturbations of $M.$

Moreover, the $\Sk_6$-simplicity implies that the pencil $j_0^1M$ of $\Sk_6$ matrices must be {\em non-singular}, that is, should contain a non-degenerate matrix: in the most generic case of this not taking place, the $\Pf \circ M$ singularity is  of the class $X_9,$ and we can again change the modulus by slightly perturbing $M.$ 

By \cite{T}, the classification of non-singular pencils of $\Sk_6$ matrices is the trivial suspension of that of $\Sq_3$ matrices:
\begin{equation} \label{Ej1}
\hspace{-12pt}
a)\!
\left( 
\begin{array} {ccc}
x & 0 & 0 \\
0 & x+y & 0 \\
0 & 0 & y 
\end{array}
\right) \ 
b)\!
\left( 
\begin{array} {ccc}
x & y & 0 \\
0 & x & 0 \\
0 & 0 & y 
\end{array}
\right) \ 
c)\! 
\left( 
\begin{array} {ccc}
x & 0 & 0 \\
0 & x & 0 \\
0 & 0 & y 
\end{array}
\right) \  
d)\!
\left( 
\begin{array} {ccc}
x & y & 0 \\
0 & x & y \\
0 & 0 & x 
\end{array}
\right) \ 
e)
\left( 
\begin{array} {ccc}
x & 0 & y \\
0 & x & 0 \\
0 & 0 & x 
\end{array}
\right)\!.
\end{equation}
The adjacency hierarchy here is 
\vspace{-23pt}

$$
\hspace{60pt} \begin{array}{ccccccc}
a & \leftarrow & b & \leftarrow & c \\
   &               &   & \nwarrow &  & \nwarrow \\
   &               &   &               & d & \leftarrow & e \   .
\end{array}
$$ The notation we are using means that, for example, the orbit of $d$ is in the closure of the orbit of $b.$

The list (\ref{Ej1}) has already appeared in \cite{BT} during the classification of simple maps $\check M: (\CC^2,0) \to (\Sq_3,0)$ in the analogous situation, 
when it was observed that the functions $\det \circ \check M$ must avoid (adjacencies to) singularities $J_{10}$ and $X_9$ (cf. Table 3).

\medskip
We now check that each of the five $j_0^1M$ options obtained by the trivial suspension of the list (\ref{Ej1}) yields as $\Sk_6$-simple singularities, the trivial suspensions of the corresponding $\Sq_3$-simple maps and only them.  The initial reductions in the cases $b)$--$e)$ below may be done by consideration of the action of the differential $d_0$ of the analogue of Arnold's spectral sequence \cite{AGV1} for which the corresponding linear part is taken as the principal part of a map.

\smallskip
$a)$ It is an easy exercise to check that in this case our matrix family $M: (\CC^2,0) \to (\Sk_6,0)$ is $\Sk$-equivalent to its linear part. 

\smallskip
$b),d)$ It is similarly easy to check that these linear parts allow us to bring our family $M$ 
to the form of the trivial suspension of a map $\hat M = (\CC^2,0) \to (\Sq_3,0)$ with its linear part staying respectively $b)$ or $d)$. Proposition \ref{Ptriv} now reduces the simple classification of the maps $M$ to that of the maps $\check M.$

\smallskip
$c)$ Some work allows to notice that
$$
j^1_0 M \sim 
\left( \begin{array} {cccccc}
0 & 0 & 0 & x & 0 & 0 \\
0 & 0 & x & 0 & 0 & 0 \\
0 & -x & 0 & 0 & 0 & 0 \\
-x & 0 & 0 & 0 & 0 & 0 \\
0 & 0 & 0 & 0 & 0 & y \\
0 & 0 & 0 & 0 & -y & 0 
\end{array}
\right) \  {\rm implies\, } \  
M \sim 
\left( \begin{array} {cccccc}
0 & \a & \b & x & 0 & 0 \\
-\a & 0 & x+\g & \d & 0 & 0 \\
-\b & -x-\g & 0 & \e & 0 & 0 \\
-x & -\d & -\e & 0 & 0 & 0 \\
0 & 0 & 0 & 0 & 0 & y \\
0 & 0 & 0 & 0 & -y & 0 
\end{array}
\right),
$$
where the Greek letters denote elements of $\frak m^2_y.$ We recognise the top left $4 \times 4$ corner here as a particular case of (\ref{E4x4}). Hence, by Section \ref{S2to4}, this corner reduces to 
$$
\left( \begin{array} {cccc}
0 & 0 & y^p & x  \\
0 & 0 & x & y^q \\
-y^p & -x & 0 & 0  \\
-x & -y^q & 0 & 0
\end{array}
\right)\!, \  2\le p \le q \ \  \Longrightarrow \ \ 
M \sim
\left( 
\begin{array} {cc}
0 & \hat M \\
-\hat M^\top & 0 
\end{array}
\right) {\rm \  where \  } \hat M=\left( \begin{array} {ccc}
y^k & x  & 0 \\
 x & y^\ell & 0 \\
 0 & 0 & y 
\end{array}
\right).
$$

$e)$ Consideration of the action of the $\Sk$-equivalence on the linear part 
allows us to reduce matrix families in this case to the form
\begin{equation}
\label{Elin-e}
M =
\left( \begin{array} {cccccc}
0 & 0 & 0 & x+\a & 0 & y \\
0 & 0 & \b & \g & x & 0 \\
0 &-\b & 0 & \d & \e & x + \eta \\
-x-\a & -\g & -\d & 0 & \zeta & 0 \\
0 & -x & -\e & -\zeta & 0 & 0 \\
-y & 0 & -x-\eta & 0 & 0 & 0 
\end{array}
\right), \quad \a, \dots, \zeta \in \frak m^2_y.
\end{equation}
Its Pfaffian avoids having a non-simple function singularity just in two cases:
\begin{itemize}
\item[(i)]
if $\d''(0) \ne 0,$ which gives $Pf \circ M \in E_7;$ 
\item[(ii)] if $\d''(0) = 0$ but 
$\left| \begin{array} {cc}
\b''(0) & \e''(0) \\
\g''(0) & \zeta''(0)
\end{array}
\right| \ne 0,$ implying $Pf \circ M \in E_8.$ 
\end{itemize}

For (i), one can show that the principal quasi-homogeneous part of $M$ (not the same as its linear part earlier) allows us to eliminate all terms of higher quasi-degrees, including the whole of $\b$ and $\zeta.$ For (ii) nearly no calculations are needed: the $\b$ and $\zeta$ may already be eliminated using the general form of $M$ by applying multiples of the operations $T_{25}$ and $T_{52},$ after which the final normal form follows from the $\Sq_3$-classification. Thus, in the $e$) case, $M$ is the trivial suspension of 
$${\rm either \quad (i) \  } \left( \begin{array} {ccc}
 x  & 0 & y \\
 0 & x &  0 \\
 y^2 & 0 & x 
\end{array}
\right) \quad  {\rm or \quad (ii) \   }
\left( \begin{array} {ccc}
x  & 0 & y \\
y^2 & x &  0 \\
0 & y^2  & x 
\end{array}
\right).
$$

This completes our classification of $\Sk_{2k}$-simple matrix families with isolated Pfaffian singularities.

\subsubsection{Maps of rank 13}     
The classification of the linear parts of such maps to $(\Sk_6,0)$ is dual in the sense of Proposition \ref{Pdual} to the classification of the rank 2 linear parts of maps considered in the previous subsection. The list (\ref{Ej1}) is the initial part of the latter, and therefore we have to start our search for rank 13 $\Sk_6$-simple matrix families by considering those maps whose linear parts may be
obtained from the linear parts (\ref{Ej1}) by taking their duals in the $\Sq_3$ sense followed by the full suspension of these duals. At some moment we will see that there are no simple rank 13 singularities with more degenerate 1-jets. 

So, for a representative of the class of rank 13 linear maps corresponding to $a)$ in (\ref{Ej1}) we can take
$$
L_a=\left( \begin{array} {cccccc}
0 & x_{12} & x_{13} & -x_{36} & x_{15} & x_{16} \\
-x_{12} & 0 & x_{23} & x_{24} &  -x_{36} & x_{26} \\
-x_{13} & -x_{23} & 0 & x_{34} &  x_{35} & x_{36} \\
x_{36} & -x_{24} & -x_{34} & 0 &  x_{45} & x_{46}  \\
-x_{15} & x_{36} & -x_{35} & -x_{45} & 0 & x_{56} \\
-x_{16} &-x_{26} & -x_{36} & -x_{46} & -x_{56} & 0 
\end{array}
\right).
$$
It is easy to check that a map germ $(\CC^{13},0) \to (\Sk_6,0)$ with this linear part reduces to just the linear part itself, has $\tau_{\Sk} =2$, and its $\Sk$-miniversal deformation is $L_a+\l_1 \E_{14} + \l_2 \E_{25},$ where the $\l_i$ are parameters of the deformation.

Therefore, any map germ with the same linear part $L_a,$ but with the source of higher dimension, reduces to the form $L_a(x)+\l_1(z) \E_{14} + \l_2(z) \E_{25},$ where $z=(z_1,\dots,z_m)$ are additional variables and $\l_1,\l_2 \in {\frak m}_z^2.$ 

We claim that moduli of the $\Sk$-classification appear here already in the quadratic parts of the $\l_i.$ 
Indeed, the problem of $\Sk$-normalisation of the pair $(j^2_0 \l_1, j^2_0 \l_2)$ while keeping the matrix $M_0$ fixed is the problem of simultaneous reduction of a pair of quadratic forms in $m$ variables by linear transformations. This problem has $m$ moduli: the set of the eigenvalues of the pair. The $L_a$ is the most general linear part in our current setting, and therefore we have

\begin{prop} \label{P13} {\em (cf. \cite{BT}, p. 757).}
There are no $\Sk_6$-simple maps from $(\CC^s,0)$ to $(\Sk_6,0)$ of rank 13 if $s>13.$
\end{prop}   

Thus, we carry on with just the classification of maps $M: (\CC^{13},0) \to (\Sk_6,0),$ and consider now those with the linear parts $L_b$ and $L_d$ constructed from the linear parts $b)$ and $d)$ from  the list (\ref{Ej1}) in the same way as $L_a$ was produced from $a).$  Straightforward calculations show that in both cases the suspending variables $x_{12}, x_{13}, x_{23}, x_{45}, x_{46}, x_{56}$ may be completely excluded from all non-linear terms of $M,$ and hence $M$ reduces to the full suspension  of a map $\check M: (\CC^{7},0) \to (\Sq_3,0).$ After that Proposition \ref{Pfull} yields a classification of $\Sk_6$-simple maps $M$ with these two types of the 1-jets.

Take now the linear family
$$
L_c=\left( \begin{array} {cccccc}
0 & x_{12} & x_{13} & -x_{25} & x_{15} & x_{16} \\
-x_{12} & 0 & x_{23} & x_{24} &  x_{25} & x_{26} \\
-x_{13} & -x_{23} & 0 & x_{34} &  x_{35} & 0 \\
x_{25} & -x_{24} & -x_{34} & 0 &  x_{45} & x_{46}  \\
-x_{15} & -x_{25} & -x_{35} & -x_{45} & 0 & x_{56} \\
-x_{16} &-x_{26} & 0 & -x_{46} & -x_{56} & 0 
\end{array}
\right).
$$
Straightforward calculations show that a matrix family $M$ with this 1-jet may be reduced to 
the form $L_c (x)+ \f(x_{12}, x_{15}, x_{24}, x_{25}, x_{45}) \E_{36}$, where $\f \in {\frak m}^2.$ The arguments of $\f$ are exactly the variables appearing in the matrix
$$
L=\left( \begin{array} {cccccc}
0 & x_{12} & -x_{25} & x_{15}  \\
-x_{12} & 0 & x_{24} &  x_{25} \\
x_{25} & -x_{24} & 0 &  x_{45} \\
-x_{15} & -x_{25} & -x_{45} & 0 
\end{array}
\right)
$$
which remains after omission of the 3rd and 6th rows and columns in $L_c,$ and these variables do not appear elsewhere in $L_c.$

We claim that $j_0^2\f$ contains moduli of the $\Sk$-equivalence.

To prove this we notice that the only available linear changes of coordinates on the domain $\CC^5$ of $\f$ may be described as follows. 
Consider the action $B^\top L B$ of elements $B \in \GL_4$ on $L.$ Take just those $B$ for which the sum of the $13$- and $24$-entries of $B^\top L B$ is zero. Such elements $B$ form a subgroup in $\GL_4$ of codimension 5, that is, of dimension 11. This subgroup acts on the $\CC^5$ by replacing the coordinates there by the corresponding entries of $B^\top L B.$ The dimension of the subgroup is smaller than that of the space of all quadratic forms in 5 variables. Hence, the $L_c$ case is indeed modal.

\medskip
Maps from $(\CC^{13},0)$ to $(\Sk_6,0)$ with the linear parts not equivalent to $L_a, L_b$ or $L_d$ are adjacent to maps with the 1-jet equivalent to $L_c.$ Therefore, there are no more simple maps from $\CC^{13}$ to $\Sk_6.$ In terms of the full suspensions, this is consistent with the $\Sq_3$-classification in part 4 of Theorem \ref{TBT}.

\begin{rem} \label{Rfence7} {\em  (cf. \cite{BT}, p. 757). 
The first non-simple principal quasi-homogeneous part of maps from $(\CC^{7},0)$ to $(\Sq_3,0)$ that appears if we try to continue the simple 3-term mini-series in part 4  of Theorem \ref{TBT} is 
$$\left( \begin{array}{ccc}
-x_{22}-x_{33} + \a x_{13}^3 & x_{12} & x_{13} \\
-x_{32} + \b x_{13}^4 & x_{22} & x_{23} \\
x_{31} & x_{32} & x_{33}
\end{array}
\right), \quad \a,\b \in \CC.
$$
One of the coefficients $\a$ and $\b$ may be normalised here, but not both simultaneously. 
}
\end{rem}

\subsection{Maps of target corank at most 1}
The corank 0 case is clear: a map germ $(\CC^s,0) \to (\Sk_{2k},0)$ of rank $N=k(2k-1)$ is a submersion, and therefore  reduces
to the normal form $\CC^N_x \times \CC^{s-N}_z \to \Sk_{2k},$ $m_{ij}=x_{ij},$ which is $\Sk$-simple and has the Tjurina number $0.$

\medskip
To establish what happens with map germs of rank $N-1,$ we first recall the $\GL_{2k}$-classification of hyperplanes in $\Sk_{2k}.$ By Proposition \ref{Pdual}, it is dual to the similar classification of lines. The latter consists of
$k$ classes, each represented by the line spanned by $\E_{1,2}+\E_{3,4} + \dots + \E_{2t-1,2t}$
for some $t=1,\dots, k.$ 
This gives us $k$ classes of hyperplanes represented by the hyperplanes 
$m_{1,2}+m_{3,4} + \dots + m_{2t-1,2t} = 0.$ Respectively, the 1-jets of rank $N-1$ map germs have normal forms
$$
-\sum\limits_{\ell=2}^t x_{2\ell-1,2\ell}\,\E_{12} 
+ \sum\limits_{1 \le i < j \le 2k, \   (i,j) \ne (1,2)}  x_{ij} \E_{ij} 
$$
It is not so difficult to show (cf. \cite{B}, p. 348)  that all non-linear terms of a map
$M: (\CC^{N-1}_x \times \CC^{s-N+1}_z,0) \to (\Sk_{2k},0)$ may in this case 
 be reduced to the form $h(X,z)\E_{12}$ where the matrix  $X$ is the lower right $2(n-t) \times 2(n-t)$ corner of $M,$ that is, $h$ is actually a function of $z$ and of
all the coordinates $x_{ij}$ in $\CC^{N-1}_x$ with $i>2t.$
In the most generic situation the restriction $h|_{X=0}$ is a Morse function, in which case a change of the $z$ coordinates makes $h(X,z) = \sum\limits_{\ell = 1}^{s-N+1} z_\ell^2 + g(X).$

We claim that the space of quadratic parts of the functions $g$ brings moduli to our $\Sk$-classification 
if $k-t>1.$ Indeed, to normalise $j^2_0g$ we need linear changes of the coordinates $X.$ Keeping $j^1_0M$ without any alternations, these changes may come only from the action $A^\top M A$ on $M$ by matrices $A = I_{2t} \oplus B$ where $B \in \GL_{2(k-t)}$: their effect on the coordinates $X$ is $B^\top X B,$ and the entries of the last matrix may be taken for new coordinates $X$. Thus we have an action of $\GL_{2(k-t)}$ on the space of the 2-forms $j^2_0 g$ in $(k-t)(2(k-t)-1)$ variables $X.$ This action has moduli when the dimension of the group is smaller than that of the space, and this is exactly when $k-t>1.$

Finally, if $t=k$ then we have no coordinates $X,$ and our matrix equivalence provides $\K$-equivalence of the functions $h$ on $\CC^{s-N+1}_z.$ This gives $\Sk_{2k}$-simple families with $h \in A_\mu, D_\mu, E_6, E_7, E_8.$  For $t=k-1,$ $X=x_{2k-1,2k}$ and it is not difficult to show that $\Sk_{2k}$-simple matrix families are classified here by the functions $h \in B_\mu, C_\mu, F_4$ on the $\CC^{s-N+2}_{x_{2k-1,2k},z}$ with the boundary $x_{2k-1,2k}=0.$

All $\Sk_{2k}$-simple maps  obtained in this subsection are clearly $\Sk$-equivalent to the full suspensions of all $\Sq_k$-simple singularities from parts (2) and (3) of Theorem \ref{TBT}.

\subsection{$\Sk$-simplicity of suspended $\Sq$-simple singularities}\label{Sskcheck}
In this section we are proving a skew-symmetric analogue of Proposition \ref{Pmatc}:

\begin{prop} \label{Pskc}
A map germ $(\CC^s,0) \to (\Sk_{2k},0)$ is $\Sk$-simple if and only if it is $\Sk_{2k}$-simple.
\end{prop}

This will complete our proof of Theorem \ref{Tmain}.

\medskip
Proving the proposition we are considering three cases:
\begin{itemize}
\item[a)] the Pfaffian of a matrix family has an isolated singularity, that is,  $s \le 6$;
\item[b)] map germs of target corank 1;
\item[c)] map germs $(\CC^{13},0) \to (\Sq_{6},0).$
\end{itemize}
In each case we will be using the suspension relations we have established with the $\Sq$-simple singularities.

As above, we call a set of non-simple equivalence classes {\em fencing} if any other non-simple singularity is adjacent to at least one element of this set. However, we are not assuming now any kind of minimality of such a set.    

\medskip
{\em Proof.} a) In this case, all $\Sk_{2k}$-simple singularities are suspensions of $\Sq_k$-simple singularities whose determinantal function has an isolated singularity, that is, $\cs \le 4.$ Analysis of the non-simplicity cutoffs in the $\Sq_k$-classification trees in \cite{BT} shows that 
in each of these classifications one can take for a fencing set of singularities the set of all map germs $\cM$ such $\det \circ \cM \in P_8, X_9, J_{10}.$ Under all three types of suspensions we have used within Table 4, the stable $\K$-equivalence classes of $\det \circ \cM$ and $\Pf \circ M$ coincide.

During the reduction of the $\Sk_{2k}$-simple classifications to the suspensions of the $\Sq_k$-simple
ones carried out earlier in Section \ref{Sred}, we also made non-simplicity cutoffs. Each of them was once again to avoid $\Pf \circ M$ having or being adjacent to the fencing function singularities.

Thus, for $s \le 6$ and $k \ge 1$, each $\Sk_{2k}$-fencing set is formed by all equivalence classes of maps whose Pfaffian is a fencing function singularity. 
Therefore, $\Sk_{2k}$-simple map germs cannot have any adjacencies to non-simple germs of lower matrix coranks. 

\medskip
In fact, the same reason is behind the validity of  Proposition \ref{Pmatc} when the determinantal function has an isolated singularity. We have already clarified that above in the $\Sq$-setting, and it may also be checked that $\Sym_k$-fencing sets in the classifications of symmetric matrix families in at most three variables done in \cite{B} are again formed by classes whose determinant is a $\K$-fencing function.

To formulate a general property,
we shall call the number of variables of a matrix family {\em low\/} if this number is at most the codimension in the corresponding matrix space of the singular locus of the variety of all degenerate matrices.
Then we have

\begin{corollary}
Assume the number of variables of a matrix family is low.
Then, in all three $\Mat$-classifications,
a matrix family is $\Mat$-simple if and only if its determinantal or respectively Pfaffian function is $\K$-simple.
\end{corollary}

b) In \cite{Gvc}, the full suspensions of the $\Sq_k$-simple singularities from part (3) of Theorem \ref{TBT} were denoted $X_\mu^{\Sk_{2k}},$ where $X_\mu$ is the corresponding $\R$- or $\R_\p$-simple (equivalently, $\K$- or $\K_\p$-simple) function germ. It follows from \cite{Gvc} (cf. Example 4.3 there) that they are adjacent to singularities of a lower matrix corank only if $X_\mu = B_\mu,C_\mu,F_4$ in which case all the adjacencies are to the singularities $Y_\nu^{\Sk_{2k-2}}$ where $Y_\nu$ is an $\R$-simple function to which $X_\mu$ deforms outside the boundary, that is, $Y_\nu = A_\nu.$

\bigskip
c) All $\Sk_6$-simple map germs of $\CC^{13}$ are full suspensions of $\Sq_3$-simple map germs of $\CC^7.$ Respectively, miniversal deformations of the former are simultaneous full suspensions of miniversal deformations of the latter (see Section \ref{Smain}). Therefore, all the $\Sk$ bifurcations may be expressed in terms of the $\Sq$ bifurcations. On the other hand, the classification of simple maps from $(\CC^7,0)$ to $(\Sq_3,0)$ was shown in \cite{Gvc} to be equivalent to the classification of simple odd functions on the plane with respect to the analogue of the usual $\K$-equivalence of functions. Here is the list of simple odd functions indexed by their Milnor numbers:
\begin{equation} \label{Eodd}
D_{2(p+1)}: \  x^2 y + y^{2p+1}, \  p \ge 1 \  \  \qquad E_8: \  x^3 + y^5 \  \   \qquad J_{10}: \  x^3 + xy^4 \  \  \qquad E_{12}: \  x^3 +y^7.  \qquad
\end{equation}
$\ZZ_2$-equivariant bifurcations of these functions translate to our current language of bifurcations of the skew-symmetric matrix families as follows:
\begin{itemize}
\item a critical point at $0 \in \CC^2_{x,y}$ corresponds to an $\Sk_6$-simple map germ of $(\CC^{13},0)$;
\item a pair of centrally symmetric critical points of type $X_\mu$ corresponds to the matrix singularity $X_\mu^{\Sk_4}.$
\end{itemize}
In the list (\ref{Eodd}), only the last two functions are not $\K$-simple in the non-equivariant sense. However, their Milnor numbers are too small to allow perturbations with two non-simple critical points.  \hfill{$\Box$}

\ 

\section{Suspending symmetric families to arbitrary square} \label{Ssymtosq} 
The similarities between the lists of $\Sym$- and $\Sq$-simple singularities obtained in respectively \cite{B} and \cite{BT} 
indicate two types of $\Sym$-to-$\Sq$ suspensions which relate these two classifications. To introduce these suspensions, consider a map germ $\hat M : (\CC^{\hat s},0) \to (\Sym_k,0),$ the standard embedding $i_k: \Sym_k \to \Sq_k$ and the standard direct sum representation $\Sq_k = \Sym_k \oplus \Sk_k.$ 

\begin{definition} \label{Dsymsq} {\em
The {\em trivial suspension} of a map $\hat M$ is
the map germ $i_k \circ \hat M : (\CC^{\hat s},0) \to (\Sq_k,0).$
The {\em full suspension} of $\hat M$ is
the map germ $\hat M \oplus V : (\CC^{\hat s},0) \oplus (\CC^{k(k-1)/2},0)\to (\Sym_k,0) \oplus (\Sk_k,0)=(\Sq_k,0),$ where $V$ is a diffeomorphism germ.}
\end{definition}

These two constructions embed the $\Sym$-simple classification of \cite{B} into the $\Sq$-simple of \cite{BT} quoted in Section \ref{Ssq}. The only $\Sq$-simple singularities which cannot be obtained in these ways are the transposed pairs of 2-variable $\Sq_2$ families in Table 1 of Theorem \ref{TBT} singled out by the braces.   Parts (2, 3, 4, 6) of that theorem are provided by the full suspensions, and the rest by the trivial ones. Some normal forms given in the theorem must be changed to make these relations obvious.   

\begin{example} {\em
The multiplicity of the determinant of  an  $\Sq$-simple map germ from $(\CC^3,0)$ to $(\Sq_2,0)$ must be 2. Therefore, such a map may be written as a two-parameter deformation induced from an $\Sq$-miniversal deformation of
the one-variable family $\left( \begin{array} {cc}
0 & z \\
-z & 0
\end{array} \right),$ for example as
 $\left( \begin{array} {cc}
\l_1& \l_2 \\
\l_2 & \l_3 
\end{array} \right)
+ \left( \begin{array} {cc}
0 & z \\
-z & 0
\end{array} \right), \  \l_1,\l_2,\l_3 \in {\frak m}_2
$.
}
\end{example}

In general, the reduction of the $\Sq$-equivalence to the $\Sym$-equivalence in the full suspension situation comes from the following.
Assume the $\Sq_k$-equivalence has brought a family of arbitrary square matrices of the matrix corank $k$ to the form $\check M = \hat M \oplus V,$ with $\hat M$ and $V$ as in Definition \ref{Dsymsq}: all upper triangle entries of the skew-symmetric matrix $V$ are $k(k-1)/2$ of the source coordinates, and the symmetric part $\hat M$ does not depend on these entries. An $\Sq$-equivalence move
$$\check M \to A^\top \check M B = A^\top \hat M B + A^\top V B$$
 (see (\ref{Esqeq})) transforms $\hat M$ and $V$ to respectively symmetric and skew-symmetric families if $A=B.$ If we do not want entries of $V$ to be involved in the new symmetric part  then the invertible family $A$ should not depend on them. After the move, the upper triangle entries of $A^\top V B$ replace those of $V$ as some of the source coordinates. 

Thus, similar to what we had in Section \ref{Ssus} for the $\Sq$-to-$\Sk$ full suspensions, we now have

\begin{prop} \label{Psymsq}
The $\Sq$-classification of those maps $\check M: (\CC^{\hat s+k(k-1)/2},0) \to (\Sq_k,0)$  which are $\Sq$-equivalent to full suspensions of maps $\hat M: (\CC^{\hat s},0) \to (\Sym_k,0)$ coincides with the $\Sym$-classification of the maps $\hat M.$
\end{prop}

Thus, one  may say that --- modulo our $\Sym$-to-$\Sq$ and $\Sq$-to-$\Sk$ suspensions --- all three simple classifications considered in this paper fit together in one injective-surjective sequence
$$
\Sym \hookrightarrow \Sq \twoheadrightarrow \Sk,
$$
with bijectivity failing exactly at the pairs of 2-parameter families of arbitrary $2 \times 2$ matrices which do not have symmetric $\Sq$-normal forms while corresponding to $\Sk$-equivalent skew-symmetric families. In particular, the $\Sym$-simple list is embedded modulo the suspensions into the $\Sk$-simple list.

\

\ 

\noindent
James William Bruce\\
Department of Mathematical Sciences\\
The University of Liverpool\\
Liverpool L69 3BX, UK\\
{\tt  billbrucesingular@gmail.com}\\

\noindent
Victor Goryunov\\
Department of Mathematical Sciences\\
The University of Liverpool\\
Liverpool L69 3BX, UK\\
{\tt goryunov@liverpool.ac.uk}\\

\noindent
Gareth Jon Haslinger \\
Ainsdale, UK \\
{\tt garethjhas@aol.com}

\end{document}